\documentstyle[12pt,twoside,amsfonts,leqno,amscd]{amsart}
\def\e{\varepsilon}
\def\l{\lambda}
\def\d{\delta}
\def\s{\sigma}
\def\f{\varphi}
\def\C{{\Bbb C}}

\def\N{{\Bbb N}}
\def\o{\omega}
\def\a{\alpha}
\def\P{{\Bbb P}}
\newcommand{\ra}{\rightarrow}
\newtheorem{thr}{Theorem}[section]
\newtheorem{THR}{Theorem}
\newtheorem{lem}[thr]{Lemma}
\newtheorem{prop}[thr]{Proposition}
\newtheorem{cor}[thr]{Corollary}
\newtheorem{rem}[thr]{Remark}
\newtheorem{defi}[thr]{Definition}
\newtheorem{exa}[thr]{Example}

\begin{document}

\title{Invariant currents and dynamical Lelong numbers}
\author{Dan Coman \ \& \ Vincent Guedj}
\thanks{D. Coman was supported by the NSF grant DMS 0140627}
\address{D. Coman: dcoman@@syr.edu, Department of Mathematics, Syracuse
University, Syracuse, NY 13244-1150, USA}
\address{V. Guedj: guedj@@picard.ups-tlse.fr, Laboratoire Emile
Picard, UMR 5580, Universit\'e Paul Sabatier, 31062 Toulouse
C\'edex 04, FRANCE}
\maketitle
\pagestyle{myheadings}
\markboth{Dan Coman \& Vincent Guedj}{Invariant currents and
dynamical Lelong numbers}

\begin{abstract}
\noindent Let $f$ be a polynomial automorphism of $\C^k$ of degree
$\l$, whose rational extension to $\P^k$ maps the hyperplane at
infinity to a single point. Given any positive closed current $S$
on $\P^k$ of bidegree (1,1), we show that the sequence
$\l^{-n}(f^n)^* S$ converges in the sense of currents on $\P^k$ to
a linear combination of the Green current $T_+$ of $f$ and the
current of integration along the hyperplane at infinity. We give
an interpretation of the coefficients in terms of generalized
Lelong numbers with respect to an invariant dynamical current for
$f^{-1}$.
\end{abstract}

\section*{Introduction}
Let $f=(P_1,\ldots,P_k):\C^k \rightarrow \C^k$ be a polynomial
automorphism of first algebraic degree $\l=\max \deg P_j \geq 2$.
We still denote by $f: \P^k \rightarrow \P^k$ the meromorphic
extension of $f$ to the complex projective space $\P^k=\C^k \cup
(t=0)$, where $(t=0)$ denotes the hyperplane at infinity.

\par The mapping $f:\P^k \rightarrow \P^k$ is not well defined on
the indeterminacy locus $I^+$, which is an algebraic subset of
$(t=0)$ of dimension $\leq k-2$. Set $X^+=f((t=0) \setminus I^+)$.
We assume throughout the paper that $X^+$ is reduced to a point
which does not belong to $I^+$. In particular $f$ is {\it weakly
regular} (see [GS]) hence it is {\it algebraically stable}: the
sequence $\l^{-n} (f^n)^*\omega$ converges in the weak sense of
currents to a positive closed current $T_+$ of bidegree $(1,1)$
such that $f^*T_+=\l T_+$ (see [S]). Here $\omega$ denotes the
Fubini-Study K\"ahler form on $\P^k$. Given $S$ a positive closed
current of bidegree $(l,l)$ on $\P^k$, we set $\|S\|:= \int_{\P^k}
S \wedge \omega^{k-l}$.

\par We assume in the sequel that $\l>\l_2(f)$, the second dynamical
degree of $f$. This allows us to construct an invariant positive
closed current $\s_-$ of bidimension $(1,1)$ which we study in
{\it section 1}. We show (Theorem 1.2) that any
quasiplurisubharmonic function is integrable with respect to the
trace measure $\sigma_- \wedge \o$. Using this we can define a
generalized Lelong number $\nu(\cdot,\s_-)$ with respect to the
dynamical weight $\s_-$ (see Definition 1.3). The dynamical
interest of these numbers lies in an invariance property
(Proposition 2.1) which we establish when $I^+$ is an
$f^{-1}$-attracting set. This last assumption has interesting
dynamical consequences (see Theorem 2.13 in [GS]).

\par Let $S$ be a positive closed current of bidegree $(1,1)$ and of
unit mass in $\P^k$. Analyzing the behavior of the bounded
sequence of currents $\l^{-n}(f^n)^*S$ is a natural problem since
it is linked with ergodic properties of the invariant current
$T_+$. This has been studied intensively in the past decade,
starting with the work of Bedford-Smillie [BS] and
 Forn\ae ss-Sibony [FS] on complex H\'enon mappings (for further
references see [S], [G1]). In the context described above, our
main result is the following:

\begin{THR} Let $f$ be a polynomial automorphism of $\C^k$ such
that $X^+$ is a point not in $I^+$. Assume that $\l>\l_2(f)$ and
that $I^+$ is an attracting set for $f^{-1}$. If $S$ is a positive
closed current on $\P^k$ of bidegree $(1,1)$ with $\|S\|=1$, then
$$\frac{1}{\l^n}\,(f^n)^* S \ra c_S [t=0]+(1-c_S)T_+,$$
in the weak sense of currents on ${\P}^k$, where
$c_S=\nu(S,\s_-)\in[0,1]$ is the generalized Lelong number of $S$
with respect to the invariant weight $\s_-$. Moreover,
$\nu(S,\s_-)>0$ if and only if the Lelong number $\nu(S,X^+)>0$.
\end{THR}

\par It should be noted that this result is new even in the case
when $f$ is a complex H\'enon mapping ($k=2$). In this case
$\s_-=T_-$ is the Green current of $f^{-1}$, hence $\nu(S,\s_-)$
is a generalized Lelong number in the sense of Demailly [D]. For
H\'enon mappings, it was shown by Bedford and Smillie that
$\l^{-n}(f^n)^*[{\mathcal C}]\rightarrow cT_+$ in $\C^2$, $c>0$,
for any algebraic curve ${\mathcal C} \subset \C^2$ (see Theorem
4.7 in [BS]). Our result can be seen as a full generalization of
this, in the sense that it yields global convergence on $\P^2$
(explaining what happens at infinity) and that it applies to any
positive closed current $S$ and in any dimension.

\par On our way to prove this theorem, we introduce an interesting
invariant probability measure $\mu_f=T_+ \wedge \s_-$ ({\it
section 1.3}). We prove Theorem 1 in {\it section 2} and we check
in {\it section 3} our hypotheses on the families of quadratic
polynomial automorphisms of $\C^3$.

\section{Invariant Lelong number}
Let $f$ be a polynomial automorphism of $\C^k$ which maps
$(t=0)\setminus I^+$ to a point $X^+\not\in I^+$ and such that $\l
> \l_2(f)$. Here $\l_2(f)$ denotes the second dynamical degree of
$f$, $\l_2(f)=\lim [\d_2(f^n)]^{1/n}$, where $\d_2(f^n)$ is the
second algebraic degree of $f^n$, i.e. the degree of $f^{-n}(L)$,
$L$ a generic linear subspace of codimension $2$ (see [S]). Under
these assumptions we can construct a positive closed current
$\sigma_-$ of bidegree $(k-1,k-1)$ and of unit mass such that
$(f^{-1})^* \sigma_-=\l \sigma_-$ (see Theorem 3.1 in [GS]).

\vskip.3cm
\subsection{Construction of $\sigma_-$}
We recall the construction of $\sigma_-$ since it is crucial for
everything that follows. Let $\Theta$ be a smooth positive closed
form of bidegree $(k-1,k-1)$ and of unit mass in $\P^k$ such that
${\rm Supp}\, \Theta \cap I^+=\emptyset$. Then ${\rm Supp}\,
(f^{-1})^* \Theta \cap (t=0)=X^+$, thus $(f^{-1})^* \Theta$ is
smooth in $\P^k \setminus \{X^+\}$. Since $(f^{-1})^* \Theta$ has
mass $\l$, there exists a current $R$ of bidegree $(k-2,k-2)$ on
$\P^k$, smooth in $\P^k\setminus\{X^+\}$, such that
$$\frac{1}{\l}\,(f^{-1})^* \Theta=\Theta+dd^c R.$$
For $W_0$ an arbitrarily small neighborhood of $X^+$ we may assume
that $0\leq R\leq C\omega^{k-2}$ in $\P^k \setminus W_0$, with a
constant $C$ depending on $W_0$. Then $0\leq(f^{-p})^* R\leq
C(f^{-p})^*\o^{k-2}$ holds in $\P^k\setminus\overline{f^p(W_0)}$.
We infer
$$\sigma^{(n)}_-:=\frac{1}{\l^n}\,(f^{-n})^*\Theta=\Theta+dd^c R_n
\longrightarrow\sigma_-:=\Theta+dd^c R_{\infty},$$ where
$R_n=\sum_{j=0}^{n-1} \l^{-j} (f^{-j})^* R$ converges to
$R_{\infty}$ in the weak sense of currents: indeed $\{R_n\}$ is an
increasing sequence of positive currents in $\P^k \setminus W_0$
(because $R \geq 0$ in $\P^k \setminus W_0$ and we can assume
$f(W_0) \subset W_0$) with bounded mass as $\l > \l_2(f)$. We will
use over and over the following facts:
$$R_n \hbox{ is smooth in }  \C^k  \hbox{ and }
R_{\infty} \geq 0 \; \hbox{ in } \; \P^k \setminus W_0.$$

\begin{rem} Let $K^-\subset\C^k$ be the set of points $z$ with
bounded backward orbit $\{f^{-n}(z)\}_{n>0}$. When $I^+$ is
$f^{-1}$-attracting it was shown in [GS] that the current
$\sigma_-$ is supported in the closure (in $\P^k$) of $K^-$, which
intersects $(t=0)$ only at the point $X^+$. This was used in
particular to show that $\sigma_-$ has full mass $1$ in $\C^k$. We
will show here that $\sigma_-$ has full mass $1$ in $\C^k$ even
when $I^+$ is not $f^{-1}$-attracting. (This occurs for certain
maps in the classes 4 and 5 from Theorem 3.1.)
\end{rem}

\par Let us recall that a function is quasiplurisubharmonic (qpsh)
if it is locally given as the sum of a plurisubharmonic function
and a smooth function.

\begin{thr}
Any quasiplurisubharmonic function is in
$L^1(\sigma_- \wedge \o)$. In particular $\sigma_-$ does not charge the
hyperplane at infinity.
\end{thr}

\begin{pf} Let $\f$ be a qpsh function and let $\f_{\e}$
be a smooth regularization of $\f$. Without loss of generality we
can assume $\f,\f_{\e} \leq 0$ and $dd^c \f, dd^c \f_{\e} \geq
-\omega$. Let $\beta$ be a smooth positive closed form of bidegree
$(1,1)$ on ${\Bbb P}^k$ vanishing in $W_0$ such that
$\o=\beta+dd^c \chi$ with $\chi \geq 0$ on ${\Bbb P}^k$. By Stokes
theorem, we have
\begin{eqnarray*}
\int (-\f_{\e}) \sigma_- \wedge \o  &=& \int (-\f_{\e}) \sigma_-
\wedge \beta+
\int (-\f_{\e}) \sigma_- \wedge dd^c \chi \\
&=& \int (-\f_{\e}) \Theta \wedge \beta+
\int dd^c (-\f_{\e}) \wedge R_{\infty} \wedge \beta
+ \int dd^c (-\f_{\e}) \wedge \chi \sigma_-\\
& \leq & \int (-\f_{\e}) \Theta \wedge \beta+
\int \omega \wedge R_{\infty} \wedge \beta
+ \int \omega \wedge \chi \sigma_-,\\
\end{eqnarray*}
since $R_{\infty} \wedge \beta \geq 0$, $\chi \sigma_- \geq 0$ and
$-dd^c \f_{\e} \leq \omega$ in $\P^k$. Letting $\e \rightarrow 0$
we get
$$0\leq\int(-\f)\sigma_-\wedge\o\leq\int(-\f)\Theta\wedge\beta+
\int\omega\wedge R_{\infty}\wedge\beta+\int \omega\wedge\chi\s_-
<+\infty,$$ since $\f$ is integrable with respect to any smooth
probability measure. In particular, when
$\f=\log|t|-\log\|[z:t]\|$ is a potential of the current of
integration along the hyperplane at infinity, this shows that the
trace measure $\sigma_- \wedge \o$ puts no mass on $(t=0)$, hence
$\sigma_-$ has full mass in $\C^k$.
\end{pf}

\vskip.3cm
\subsection{Dynamical Lelong number}
Let $S$ be a positive closed current of bidegree $(1,1)$ and unit
mass on ${\Bbb P}^k$, so $S=\omega+dd^c \f$ for some qpsh function
$\f$. It follows from Theorem 1.2 that the probability measure
$S\wedge\sigma_-:=\o\wedge\sigma_-+dd^c(\f \, \sigma_-)$ is well
defined.

\begin{defi}
The generalized Lelong number of $S$ with respect to the invariant
current $\s_-$ is $\nu(S,\s_-):=S \wedge \s_-(\{X^+\})$.
\end{defi}

\par The following convergence result will help to compute
generalized Lelong numbers.

\begin{thr}
Let $S$ be a positive closed current of bidegree $(1,1)$ on $\P^k$.
Then
$$S \wedge \sigma_-^{(n)} \rightarrow S \wedge \sigma_-,$$
in the weak sense of measures on $\P^k$.
\end{thr}

\begin{pf} We can assume $S$ has mass $1$, hence
$S=\o+dd^c \f$, where $\f \leq 0$ is qpsh. We are going to show
that $\f \sigma_-^{(n)} \rightarrow \f \sigma_-$ in $\P^k
\setminus X^+$.

\par Observe first that the currents $\f \sigma_-^{(n)}$ have
uniformly bounded mass in $\P^k$: arguing as in the proof of
Theorem 1.2, we get
$$0\leq\int(-\f)\sigma_-^{(n)}\wedge\o
\leq\int(-\f)\Theta\wedge\beta+\int\o\wedge R_n\wedge\beta+\int\o
\wedge\chi\sigma_-^{(n)}\leq C<+\infty$$ since $R_n$ increases to
$R_{\infty}$ in $\P^k \setminus W_0$ and $\sigma_-^{(n)}$ has
bounded total mass.

\par Let $\nu$ be a cluster point of $\{\f \sigma_-^{(n)}\}$. Let
$\{\f_{\e}\}$ be a sequence of smooth qpsh functions decreasing
pointwise to $\f$. Then $\f \sigma_-^{(n)} \leq \f_{\e}
\sigma_-^{(n)}$, hence $\nu \leq \f_{\e} \sigma_-$. Letting $\e
\rightarrow 0$ yields $\nu \leq \f \sigma_-$. To get equality, it
suffices to show that the total mass of $(-\f) \sigma_-$ dominates
that of $-\nu$. Recall that $\sigma_-^{(n)}=\Theta+dd^c R_n$,
where $R_n=\sum_{j=0}^{n-1} \l^{-j} (f^{-j})^* R$, and $R$ is
smooth in $\P^k \setminus\{X^+\}$. Up to now, we have chosen $R
\geq 0$ in $\P^k \setminus W_0$. Here it is actually more
convenient to choose a negative potential. Set $T=R-C \o^{k-2}$,
where $C$ is a positive constant so large that $T \leq 0$ in $\P^k
\setminus W_0$. Then $\sigma_-^{(n)}=\Theta+dd^c T_n$, where
$T_n=\sum_{j=0}^{n-1} \l^{-j} (f^{-j})^* T$ is a sequence of
negative currents in $\P^k \setminus W_0$ decreasing to
$T_{\infty}$. Set
$$\hat{T}_n:=\sum_{j\geq n}\frac{1}{\l^j}\,(f^{-j})^*T\leq0
\hbox{ in }\P^k\setminus W_0,$$ so that $\sigma_-
-\sigma_-^{(n)}=dd^c \hat{T}_n$. Let $\beta$ be a smooth closed
form of bidegree $(1,1)$ on $\P^k$ vanishing in $W_0$ and strictly
positive in $\P^k \setminus \overline{W_0}$. Using $-\hat{T}_n
\wedge \beta \geq 0$ in $\P^k$, we get
\begin{eqnarray*}
\int (-\f_\e) \sigma_- \wedge \beta &=& \int (-\f_\e)
\sigma^{(n)}_- \wedge \beta+
\int (-\f_\e) dd^c \hat{T}_n \wedge \beta \\
&=& \int (-\f_\e) \sigma_-^{(n)} \wedge \beta+
\int  dd^c \f_\e \wedge (-\hat{T}_n ) \wedge \beta \\
&\geq&  \int (-\f_\e) \sigma_-^{(n)} \wedge \beta- \int \o \wedge
(-\hat{T}_n) \wedge \beta.
\end{eqnarray*}
As $\e\ra0$
$$\int(-\f)\sigma_-\wedge\beta\geq\int(-\f)\sigma_-^{(n)}\wedge\beta+\int\o
\wedge\hat{T}_n\wedge\beta.$$ Now $\hat{T}_n \rightarrow 0$ as
$n\ra+\infty$, hence $\int (-\f) \sigma_- \wedge \beta \geq \int
(-\nu)\wedge \beta$. This shows that $\nu=\f \sigma_-$ in $\P^k
\setminus W_0$, hence in $\P^k \setminus X^+$ since $W_0$ is an
arbitrarily small neighborhood of $X^+$.

\par It follows that $S \wedge \sigma_-^{(n)} \rightarrow S \wedge
\sigma_-$ in $\P^k \setminus X^+$. Since these are all probability
measures, we actually get $S \wedge \sigma_-^{(n)} \rightarrow S
\wedge \sigma_-$ on $\P^k$.\end{pf}

\begin{exa} If $\mu_n=\sigma_-^{(n)} \wedge [t=0]$ then
$\limsup\mu_n(\{X^+\})\leq\nu([t=0],\sigma_-)\leq 1$ by Theorem
1.4. Now $\mu_n(\{X^+\})=1$ because $\sigma_-^{(n)}$ clusters at
infinity only at $X^+$. Therefore $\nu([t=0],\sigma_-)=1$, i.e.
$[t=0] \wedge \sigma_-$ is the Dirac mass at the point $X^+$. At
the other end, observe that $T_+$ vanishes in a neighborhood of
$X^+$ which is an attracting fixed point, so
$\nu(T_+,\sigma_-)=0$.
\end{exa}

\par Regular automorphisms were introduced by Sibony [S] and studied in
[S], [GS]. These are automorphisms such that $I^+ \cap
I^-=\emptyset$. In this case $f^{-1}$ is algebraically stable, so
there is a well defined invariant Green current $T_-$ for $f^{-1}$
(see [S]).

\begin{prop}
Assume $f$ is a regular automorphism. Then $\sigma_-=T_-^{k-1}$, so
$\nu(S,\sigma_-)$ is the Demailly number of $S$ with respect
to the weight $T_-$. In this case,
$$\nu(S,\sigma_-)>0 \;\hbox{ if and only if }\;\nu(S,X^+)>0,$$
where $\nu(S,X^+)$ denotes the standard Lelong number at the point
$X^+$.
\end{prop}

\begin{pf} When $f$ is a regular automorphism as defined
in [S], the inverse $f^{-1}$ has first algebraic degree $d_-$ such
that $d_-^{k-1}=\l$ (recall that $X^+$ is a point), and
$\l_2(f)=d_-^{k-2}<\l$. Note also that in this case $I^+=X^-$ is
an $f^{-1}$-attracting set. We refer the reader to [S] for the
construction of $T_-=\o+dd^c g_-$, the Green current of bidegree
$(1,1)$ for $f^{-1}$. It follows from the extension of the
Bedford-Taylor theory of Monge-Amp\`ere operators that $T_-^{k-1}$
is well defined and equals $\lim \l^{-n} (f^{-n})^* (\o^{k-1})$
(see [D], [S]). Thus $T_-^{k-1}=\lim \l^{-n} (f^{-n})^*
\Theta=\sigma_-$ since $\Theta=\o^{k-1}+dd^c \a$, where $\a$ is a
smooth form of bidegree $(k-2,k-2)$, hence
$\|(f^{-n})^*(\a)\|=O(d_-^{n(k-2)})=o(\l^n)$. Note also that
$T_-^k$ is well defined and equals the Dirac mass at the point
$X^+=I^-$. This is a situation where the Jensen type formulas of
Demailly simplify and give a nice understanding of the generalized
Lelong numbers $\nu(S,T_-^{k-1})$.

\par The potential $g_-$ of $T_-$ is obtained as $g_-=\sum_{n \geq 0}
d_-^{-n} \phi_- \circ f^{-n}$, where $d_-^{-1}(f^{-1})^*
\omega=\omega+dd^c \phi_-$. Observe that $g_-$ has positive Lelong
number at $X^+ = I^-$, hence $g_-(z)\leq \gamma_1 \log
dist(z,X^+)+C$.

\par We also have control from below, $\gamma_2\log dist(z,X^+)-C\leq
g_-(z)$. This follows from a Lojasiewicz type inequality, since
$$\phi_-(z)=2^{-1}\log[|Q_0(z)|^2+\ldots+|Q_k(z)|^2]+
\hbox{ smooth term near } X^+,$$ where $Q_j$ are polynomials such
that $\bigcap Q_j^{-1}(0)=X^+$. It follows from the
Nullstellensatz that $|Q_0(z)|^2+\ldots+|Q_k(z)|^2 \geq
dist(z,X^+)^{\alpha}$ near $X^+$ for some exponent $\alpha>0$. As
$X^+$ is an attracting fixed point for $f$, we get $dist(f(z),X^+)
\leq c\,dist(z,X^+)$ for all $z \in \C^k$, hence
$dist(f^{-n}(z),X^+) \geq c^{-n} dist(z,X^+)$. Therefore $g_-(z)
\geq \gamma_2 \log dist(z,X^+)-C$ with $\gamma_2=2^{-1} \alpha
d_-/(d_- -1)$.

\par We conclude by the first comparison theorem of Demailly [D] that
$\nu(S,\sigma_-)>0$ if and only if $\nu(S,X^+) >0$. \end{pf}

\begin{rem}
For regular automorphisms $T_-^k$ is the Dirac mass at the point
$X^+=I^-$, thus $\nu(T_-,\s_-)=1$. It is an interesting question
to characterize the closed positive currents $S \sim \o$ such that
$\nu(S,\sigma_-)=1$.
\end{rem}

\begin{thr}
Let $S$ be a positive closed current of bidegree $(1,1)$
on $\P^k$.
\par 1) The sequence of currents $S \wedge R_n$ is well defined and
convergent in $\C^k$. Set $S\wedge R_{\infty}:=\lim S\wedge R_n$
in $\C^k$. Then
$$S\wedge\s_-=S\wedge\Theta+dd^c(S\wedge R_{\infty})\hbox{ in }\C^k.$$
\par 2) Assume $S_n \rightarrow S$, where $S_n$ are positive closed
currents of bidegree $(1,1)$ on $\P^k$. Then $ S_n \wedge \s_-
\longrightarrow S \wedge \s_- \hbox{ in } \C^k. $ Moreover when
$I^+$ is $f^{-1}$-attracting, then $ S_n \wedge \sigma_-
\longrightarrow S \wedge \sigma_-$ on $\P^k$.
\end{thr}

\begin{cor}
If $I^+$ is $f^{-1}$-attracting, the mapping  $S \mapsto
\nu(S,\sigma_-)$ is upper semicontinuous.
\end{cor}

\begin{pf} Let $S_n \ra S$. Then $S_n\wedge\sigma_-\ra
S\wedge\sigma_-$ on $\P^k$, so $\limsup S_n \wedge
\sigma_-(\{X^+\})\leq S \wedge \s_-(\{X^+\})$. \end{pf}

\begin{lem}
Let $S$ be a positive closed current of bidegree $(1,1)$ on $\P^k$
and let $\theta$ be a positive closed current of bidimension $(1,1)$
which is smooth in an open subset $\Omega$ of $\P^k$. Then
$$0\leq\int_{\Omega}S\wedge\theta\leq\|S\|\cdot\|\theta\|,$$
where $\|S\|=\int_{\P^k} S \wedge \o^{k-1}$ and
$\|\theta\|=\int_{\P^k}\theta\wedge\o$.
\end{lem}

\begin{pf} Since $\P^k$ is homogeneous (i.e. $Aut(\P^k)$
acts transitively on $\P^k$), we can regularize $S$ in the
following sense: there exist smooth positive closed currents
$S_{\e}$ of bidegree $(1,1)$ on $\P^k$ such that
$\|S_{\e}\|=\|S\|$ and $S_{\e} \rightarrow S$ on $\P^k$ (see [H]).
Therefore $S_{\e} \wedge \theta \rightarrow S \wedge \theta$ in
$\Omega$, hence
$$0\leq\int_{\Omega}S\wedge\theta\leq\liminf_{\e\ra0}
\int_{\Omega}S_{\e}\wedge\theta\leq\liminf_{\e\ra0}\int_{\P^k}S_{\e}
\wedge\theta = \|S\|\cdot\|\theta\|.$$\end{pf}

\vskip.1cm \noindent {\em Proof of Theorem 1.8.} Let $S$ be a
positive closed current of bidegree $(1,1)$ on $\P^k$. Recall that
$\sigma_-=\Theta+dd^c R_{\infty}$, where
$R_{\infty}=R_n+\hat{R}_n=\lim R_n$, $R_n=\sum_{j=0}^{n-1} \l^{-j}
(f^{-j})^* R$ being smooth in $\C^k$. Therefore $S \wedge R_n$ is
a well defined current of bidimension $(1,1)$ which is positive in
$\C^k \setminus W_0$. We estimate its mass in $\C^k \setminus
W_0$: if $S_{\e}$ is a regularization of $S$ as in the proof of
Lemma 1.10, then
\begin{eqnarray*}
 & & 0\leq\int_{\C^k \setminus\overline W_0}S\wedge R_n\wedge\o
\leq\liminf_{\e\ra0}\int_{\C^k \setminus\overline W_0}S_\e\wedge
R_n\wedge \o\leq \\
 & & C\liminf_{\e \rightarrow 0}\sum_{j=0}^{n-1} \frac{1}{\l^j}
\int_{\P^k} S_{\e} \wedge (f^{-j})^* \o^{k-2} \wedge \o \leq C
\|S\| \sum_{j \geq 0} \frac{\d_2(f^j)}{\l^j} <+\infty,
\end{eqnarray*}
where $C>0$ is a constant depending on the fixed neighborhood
$W_0$. This shows that the increasing sequence $\{S \wedge R_n\}$
is convergent in $\C^k \setminus\overline W_0$. Observe that the
sequence $\{R_n-R_{p}\}_{n\geq p}$ is positive and increasing in
$\C^k \setminus\overline{f^p(W_0)}$. Thus $S \wedge R_n$ converges
in $\C^k \setminus\overline{f^p(W_0)}$, for all $p$, hence in
$\C^k$, as $\overline{f^p(W_0)}\searrow X^+$. Set $S \wedge
R_{\infty}:=\lim S \wedge R_n$ in $\C^k$. Then
$$S\wedge\Theta+dd^c(S\wedge R_{\infty})=\lim[S\wedge
\Theta+dd^c(S \wedge R_n)]=\lim S\wedge\s_-^{(n)}=S \wedge \s_-,$$
by Theorem 1.4. This proves 1).

\par Let now $S_n,\,S$ be positive closed currents of bidegree $(1,1)$
on $\P^k$ such that $S_n \rightarrow S$. Since $R_N$ is smooth in
$\C^k$, we get $S_n \wedge R_N \rightarrow S \wedge R_N$ for all
fixed $N$. We want to show that $S_n \wedge R_{\infty} \rightarrow
S \wedge R_{\infty}$. It is sufficient to get an estimate on
$\|S_n \wedge\hat R_N\|_{\C^k \setminus\overline{f^N(W_0)}}$ which
is uniform in $n$. This is the following
$$0\leq\int_{\C^k \setminus\overline{f^N(W_0)}}S_n \wedge
\hat{R}_N \wedge\o\leq C\liminf_{\e\ra0}\sum_{j \geq N}
\frac{1}{\l^j}\int_{\P^k} S_n^{\e} \wedge (f^{-j})^* \o^{k-2}
\wedge\o\leq C'\sum_{j \geq N} \frac{\d_2(f^j)}{\l^j}\;,$$ where
$S_n^{\e}$ is a regularization of $S_n$ and the last inequality
follows from Lemma 1.10 and the fact that the sequence of norms
$\|S_n^{\e}\|=\|S_n\|$ is bounded. Therefore $S_n \wedge
R_{\infty} \rightarrow S \wedge R_{\infty}$ in $\C^k$, hence
$$S_n \wedge \s_-=S_n \wedge \Theta+dd^c(S_n \wedge R_{\infty})
\rightarrow S \wedge \s_- \;\rm{in}\;\C^k.$$

\par When $I^+$ is $f^{-1}$-attracting, the current $\s_-$ clusters
at infinity only at $X^+$. Since $S_n \wedge \s_-$ and $S \wedge
\s_-$ are positive measures on $\P^k$ supported in ${\rm
Supp}\,\s_-$ and $\|S_n\wedge\s_-\|=\|S_n\|\ra\|S\wedge\s_-\|$, we
infer in this case that $S_n \wedge \s_- \rightarrow S \wedge
\s_-$ on $\P^k$. $\;\;\;\Box$

\vskip.3cm
\subsection{Invariant measure}
In this section we introduce and study a dynamically interesting
probability measure.

\begin{defi}
We write $T_+=\o+dd^cg_+$ and set
$$\mu_f=T_+ \wedge \s_-:=\o \wedge \s_-+dd^c(g_+\s_-).$$
\end{defi}

\par Note that this measure is well defined thanks to
Theorem 1.2. It is clearly a probability measure since
$\int_{\P^k}\o\wedge\s_-=1$.

\par We have $T_+=0$ in the basin of attraction of $X^+$.
If $I^+$ is $f^{-1}$-attracting then the support of $\s_-$
intersects $(t=0)$ only at $X^+$ (see Remark 1.1). It follows that
in this case $\mu_f$ has compact support in $\C^k$ and it is
invariant, i.e. $f_*\mu_f=\mu_f$.

\par When $f$ is a regular automorphism, we have $\s_-=T_-^{k-1}$,
so $PSH(\C^k)\subset L^1(\mu_f)$, by the Chern-Levine-Nirenberg
inequalities. More generally, when there exists partial Green
functions for $f^{-1}$, one also gets $PSH(\C^k) \subset
L^1(\mu_f)$ (see section 4.2 in [GS]). This requires however
delicate estimates on the growth of $f^{-1}$ near $I^+$. We now
establish in the spirit of [G2] the following integrability
result:

\begin{thr}
If $I^+$ is $f^{-1}$-attracting and $\f$ is a
quasiplurisubharmonic function on $\P^k$, then $\f \in
L^1(\mu_f)$.
\end{thr}

\begin{pf} We can assume without loss of generality that
$\f<0$ and $dd^c\f \geq-\o$. Let $\f_\e<0$ be qpsh functions which
decrease pointwise to $\f$ such that $dd^c\f_\e\geq-\o$. The
current $T_+=\o+dd^c g_+$ has potential $g_+<0$ which is
continuous in $\P^k \setminus I^+$. Since $I^+$ is an attracting
set for $f^{-1}$, the current $\s_-$ vanishes in a neighborhood
$V_0$ of $I^+$. If $A=\|g_+\|_{L^{\infty}(\P^k \setminus V_0)}$
then $(g_+ +A) \s_-\geq 0$ on $\P^k$. We get
\begin{eqnarray*}
\int (-\f_{\e})d\mu_f & = & \int(-\f_{\e})\o\wedge\s_-+
\int(-\f_{\e})dd^c((g_++A)\s_-)\\
 & = & \int(-\f_{\e})\o\wedge\s_-+ \int
dd^c(-\f_{\e})\wedge((g_++A)\s_-)\\
 & \leq & \int(-\f_{\e})\o\wedge\s_-+\int(g_++A)\o\wedge\s_-
\leq A+\int(-\f_{\e})\o\wedge\s_-.
\end{eqnarray*}
The conclusion follows by letting $\e\ra0$ and using Theorem 1.2.
\end{pf}

\begin{rem}
If $u$ is a plurisubharmonic (psh) function defined in a
neighborhood of the support of $\mu_f$, then $|u|^{\alpha}\in
L^1(\mu_f)$ for every $\alpha\in(0,1/k)$. {\rm Indeed, by Theorem
1.12 psh functions of logarithmic growth are integrable with
respect to $\mu_f$. The claim is straightforward using the
following result of El Mir and Alexander-Taylor (see \cite{AT}):
If $u\leq-1$ is psh in a ball $B(z_0,R)\subset\C^k$ and $r<R$,
$0<\epsilon<1/k$, then there exists a psh function $v$ on $\C^k$
of logarithmic growth such that $v\leq-|u|^{1/k-\epsilon}$ on
$B(z_0,r)$.}
\end{rem}

\section{Equidistribution towards $T_+$}
The purpose of this section is to prove Theorem 1 stated in the
Introduction.

\begin{pf} The proof of the theorem is divided into four steps.

\vskip.1cm\noindent{\bf Step 1: Normalization of potentials.} By
Siu's theorem, we can write
\begin{equation}\frac{1}{\l^n}\,(f^n)^*S=c_n[t=0]+(1-c_n)S_n,
\end{equation}
where $c_n \in [0,1]$, and $S_n$ are positive closed currents of
bidegree $(1,1)$ and unit mass which do not charge $(t=0)$. Since
$f^*[t=0]=\l [t=0]$, the sequence $\{c_n\}$ is increasing. Let
$c_S$ denote its limit. If $c_S=1$ the convergence statement of
the theorem is proved, so we assume hereafter that $c_S<1$.

\par We write $S=\o+dd^c v_0$, where the potential $v_0$ is uniquely
determined up to additive constants. Using Theorem 1.12, we can
normalize it so that $\int v_0 d\mu_f=0$. Similarly, we fix
potentials $S_n=\o+dd^c v_n$, $T_+=\o+dd^c g_+$, $[t=0]=\o+dd^c
\f_{\infty}$ such that $\int v_n d\mu_f=\int g_+ d\mu_f=\int
\f_{\infty} d\mu_f=0$. If $\l^{-n}(f^n)^* \o=\o+dd^cg_+^{(n)}$,
$\int g_+^{(n)}d\mu_f=0$, then $g_+^{(n)}\ra g_+$ in $L^1(\P^k)$
and $\l^{-n}(f^n)^*\o\ra T_+$. The desired convergence follows if
we show that $\l^{-n} v_0 \circ f^n \rightarrow
c_S(\f_{\infty}-g_+)$ in $L^1(\P^k)$.

\par Pulling back (1) (with $n=p$) by $f^n$ yields
\begin{eqnarray*}
\frac{1}{\l^{n+p}}\,(f^{n+p})^*S&=&
c_p[t=0]+(1-c_p)\frac{1}{\l^n}\,(f^n)^*S_p \\
&=&c_p[t=0]+(1-c_p)\frac{1}{\l^n}\,(f^n)^*\o+
(1-c_p)dd^c\left(\frac{1}{\l^n}\,v_p \circ f^n\right). \\
\end{eqnarray*}
Using our normalization and the fact that $\mu_f$ is invariant, we infer
\begin{equation}
\frac{1}{\l^{n+p}}\,v_0 \circ f^{n+p}= c_p
(\f_{\infty}-g_+^{(n)})+(g_+^{(n)}-g_+^{(n+p)})
+(1-c_p)\frac{1}{\l^n} \,v_p \circ f^n.
\end{equation}

\vskip.1cm \noindent {\bf Step 2: Control of the Lelong numbers.}
Since $f^n$ is a biholomorphism in $\C^k$, it follows from (1)
that for all $n \in \N$ and $z \in \C^k$,
$$\nu((1-c_n)S_n,z)=\frac{1}{\l^n}\,\nu((f^n)^*S,z)=
\frac{1}{\l^n}\,\nu(S,f^n(z))\leq\frac{1}{\l^n}\;,$$ hence
$\sup_{z \in \C^k} \nu(S_n,z) \leq (1-c_S)^{-1} \l^{-n} \ra0$.

\par Pulling back $(1)$ by $f$ we get
\begin{equation}
\frac{1}{\l}\,f^*S_n=\frac{c_{n+1}-c_n}{1-c_n}\,[t=0]+
\frac{1-c_{n+1}}{1-c_n}\,S_{n+1}.\end{equation} Since $S_{n+1}$
does not charge $(t=0)$, we have for a generic point $z\in(t=0)$
$$\nu(S_n,X^+)=\nu(S_n,f(z)) \leq \nu(f^*S_n,z)=\l\,
\frac{c_{n+1}-c_n}{1-c_n} \leq\l\,\frac{c_S-c_n}{1-c_S}\;.$$

\par If $z \in (t=0) \setminus I^+$, it follows from [F]
and [K] that there is an upper estimate $\nu(f^*S_n,z) \leq
c_{f,z} \nu(S_n,f(z))$, where $z \mapsto c_{f,z}$ is locally upper
bounded. Fix $V_0$ a small neighborhood of $I^+$ and set
$C_{V_0}=\sup_{z \in (t=0) \setminus V_0} c_{f,z}$. Using (3)
again, we get for all $z \in (t=0) \setminus V_0$,
$$\frac{1-c_{n+1}}{1-c_{n}}\,\nu(S_{n+1},z) \leq
\frac{1}{\l}\,\nu(f^*S_n,z) \leq \frac{C_{V_0}}{\l}\,\nu(S_n,X^+)
\leq C_{V_0} \frac{c_S-c_n}{1-c_S}\;.$$

\par We conclude that $\sup_{z \in \P^k \setminus V_0} \nu(S_n,z)
\rightarrow 0$ as $n\ra+\infty$.

\vskip.1cm \noindent {\bf Step 3: Volume estimates.} We have to
prove that
$$w_n:=\l^{-n} v_0 \circ f^n \rightarrow c_S (\f_{\infty}-g_+).$$
Observe first that the sequence $\{w_n\}$ is relatively compact in
$L^1(\P^k)$. Indeed
$$\l^{-n}(f^n)^*S=\l^{-n}(f^n)^*\o+dd^c(w_n)=
\o+dd^c(g_+^{(n)}+w_n),$$
so $w_n+g_+^{(n)}$ are qpsh functions whose curvature is uniformly
bounded from below by $-\o$. Since $g_+^{(n)} \rightarrow g_+$ and
$w_n \leq C \l^{-n}$, the sequence $\{w_n+g_+^{(n)}\}$ is
uniformly upper bounded on $\P^k$. So either this sequence
converges uniformly to $-\infty$, or it is relatively compact in
$L^1(\P^k)$ (see Appendix in [G1]). The former cannot happen since
$\int (w_n+g_+^{(n)}) d\mu_f=0$. Thus it suffices to show that
$w_n$ converges in measure to $c_S (\f_{\infty}-g_+)$. It follows
from (2) that
\begin{eqnarray*}
 & & w_{n+p}-c_S (\f_{\infty}-g_+)= \\
 & & (c_p-c_S)(\f_{\infty}-g_+^{(n)})+c_S(g_+-g_+^{(n)})+
     (g_+^{(n)}-g_+^{(n+p)})+(1-c_p)\l^{-n}v_p \circ f^n.
\end{eqnarray*}

\par Let $\e>0$. Choose a small neighborhood $V_0$ of $I^+$
and fix $p$ so large that
$$\sup_{z \in \P^k \setminus V_0} \nu(S_p,z) \leq \e^2 \; \; {\rm
and} \; \; |c_p-c_S| \|\f_{\infty}-g_+^{(n)}\|_{L^1(\P^k)}<\e^2,
\;\forall\,n \in \N.$$ By Chebyshev's inequality ${\rm
Vol}(|(c_p-c_S)(\f_{\infty}-g_+^{(n)})|>\e/3)<3 \e$. Since
$g_+^{(n)}\ra g_+$ in $L^1(\P^k)$, we have for $n$ large ${\rm
Vol}(|c_S(g_+-g_+^{(n)})+(g_+^{(n)}-g_+^{(n+p)})|>\e/3)<\e$.
Observe that
\begin{eqnarray*}
\lefteqn{ {\rm Vol}(|w_{n+p}-c_S (\f_{\infty}-g_+)|>\e) \leq
{\rm Vol}(|(c_p-c_S)(\f_{\infty}-g_+^{(n)})|>\e/3)+} \\
& & {\rm Vol}(|c_S(g_+-g_+^{(n)})+(g_+^{(n)}-g_+^{(n+p)})|>\e/3)+
{\rm Vol}((1-c_p)|\l^{-n}v_p \circ f^n|>\e/3),
\end{eqnarray*}
Since $v_p$ is bounded above on $\P^k$, it remains to show that
$${\rm Vol}(|\l^{-n}v_p \circ f^n|>\e/3)=
{\rm Vol}(\l^{-n}v_p \circ f^n<-\e/3)<C \e,$$ for all $n$
sufficiently large.

\par Since $I^+$ is $f^{-1}$-attracting, there exist arbitrarily
small neighborhoods $V_0$ of $I^+$ such that $f(\P^k \setminus
V_0) \subset \P^k \setminus V_0$. Set
$$\Omega_n^{\e}:=\{ z \in \P^k \setminus V_0:\,\l^{-n} v_p \circ
f^n(z) <-\e/3 \}.$$ We have $f^n(\Omega_n^{\e}) \subset \{ z \in
\P^k \setminus V_0:\, v_p(z)<-\e \l^n /3\}$. It follows from [G1]
that there exists $C_1>0$ such that
$${\rm Vol}(f^n(\Omega_n^{\e})) \geq \exp \left( -\frac{C_1
\l^n}{{\rm Vol}(\Omega_n^{\e})} \right).$$ On the other hand, by
Skoda's integrability theorem (see [K]) there exists $C_{\e}>0$
such that
\begin{eqnarray*}
{\rm Vol}\left( \{ z \in \P^k \setminus V_0:\, v_p(z)<-\e \l^n/3
\} \right) & \leq & C_{\e} \exp \left( -\frac{\e \l^n}{3\sup_{z
\in\P^k \setminus V_0} \nu(S_p,z)}\right) \\
 & \leq & C_{\e}\exp\left(-\frac{\l^n}{3\e}\right).
\end{eqnarray*}
Thus ${\rm Vol}(\Omega_n^{\e})\leq4C_1 \e$ for all $n>N(\e)$.

\par We conclude that $w_n\ra c_S(\f_{\infty}-g_+)$ in measure
on $\P^k\setminus V_0$. As $V_0$ was an arbitrarily small
neighborhood of $I^+$, the convergence in measure holds on $\P^k$.

\vskip.1cm \noindent {\bf Step 4: Interpretation of $c_S$.} We
have shown that $\l^{-n}(f^n)^*S \rightarrow c_S[t=0]+(1-c_S)T_+$.
It follows from [G1] that $c_S>0$ if and only if $\nu(S,X^+)>0$.
Assume now that $I^+$ is $f^{-1}$-attracting. We show below
(Proposition 2.1) that $\nu((f^n)^*S,\sigma_-)=\l^n
\nu(S,\sigma_-)$. It then follows from Example 1.5 that
$$\nu(S,\s_-)=\nu(\l^{-n}(f^n)^*S,\s_-)=c_n+\nu(\hat{S}_n,\s_-),$$
where $\hat{S}_n=(1-c_n)S_n \rightarrow (1-c_S)T_+$. Since
$\nu(T_+,\sigma_-)=0$, we infer from the upper semicontinuity
property (Corollary 1.9) that $\nu(\hat{S}_n,\sigma_-) \rightarrow
0$, hence $c_S=\nu(S,\sigma_-)$. \end{pf}

\begin{prop}{\bf(Transformation rule)}
$\nu(f^*S,\sigma_-)=\l \nu(S,\sigma_-)$.
\end{prop}

\begin{pf} Let $S_j$ be a sequence of smooth closed
positive currents of bidegree (1,1) with smooth potentials which
decrease pointwise to a potential of $S$. Let $W$ be a small
neighborhood of $X^+$ so that $f(W)\subset\subset W$. Note that
$f(W)=f(W\cap\C^k)\cup X^+$. Since $f^*S_j$ is smooth in $W$ and
$\s_-$ does not charge $(t=0)$ (Theorem 1.2) we have
$$\int_W f^* S_j\wedge \sigma_-=\int_{W\cap\C^k}f^*
S_j\wedge\sigma_-= \int_{f(W)\cap{\Bbb C}^k}S_j\wedge(f^{-1})^*
\sigma_- =\l \int_{f(W)} S_j \wedge \sigma_-.$$ By the monotone
convergence theorem, one has $S_j\wedge\s_-\ra S\wedge\s_-$ and
$f^*S_j\wedge\s_-\ra f^*S\wedge\s_-$. We infer $\int_W f^*S\wedge
\sigma_-\leq\l\int_{\overline{W}}S\wedge\sigma_-$, hence
$\nu(f^*S,\s_-)\leq\l\nu(S,\s_-)$.

\par For the opposite inequality, observe that the restriction of
$f^{-1}:K^-\ra K^-$ extends continuously at infinity by setting
$f^{-1}(X^+)=X^+$. This shows $f$ is an open mapping on $\overline
{K^-}$, so there is a ball $B \subset W$ centered at $X^+$ such
that $\overline{K^-}\cap B \subset f(W)$. Therefore $ \int_{W} f^*
S_j \wedge \sigma_- \geq \l \int_B S_j \wedge \sigma_-$, which
yields
$$\int_{\overline W}f^*S\wedge\sigma_-\geq \l\int_B
S\wedge\sigma_-\geq\l\nu(S,\sigma_-).$$ The desired inequality
follows by shrinking $W\searrow X^+$.
\end{pf}

\begin{rem}
We showed in the proof of Theorem 1 that if $S=\o+dd^cv_0$ then
$\l^{-n}v_0\circ f^n\ra c_S(\phi_\infty-g_+)$ in $L^1(\P^k)$. Let
$G^+(z,t)$, $(z,t)\in\C^{k+1}$, be the logarithmically homogeneous
Green function of $f$. The function $h[z:t]=\log|t|-G^+(z,t)$ is
well defined on $\P^k$ and $h=\phi_\infty-g_++c$ for some constant
$c$. Since $h\circ f=\l h$ and $f_*\mu_f=\mu_f$ we have $\int
h\,d\mu_f=0$, so $\phi_\infty-g_+=h$.
\end{rem}

\begin{rem}
The convergence $\l^{-n}(f^n)^* S\ra c_S[t=0]+(1-c_S)T_+$ holds
without the hypotheses $\l>\l_2(f)$ and $I^+$ is
$f^{-1}$-attracting. A proof can be given in the basin of $X^+$ by
a similar argument, and on the complement of this basin one can
conclude as in the proof of Theorem 2.7 in [G1]. However in this
case we do not have an interpretation for $c_S$. As an example,
our convergence theorem holds for the maps $f$ and $f^{-1}$, where
$f(x,y,z)=(P(y)+az,Q(y)+bx,y)$, $\deg(P)=\deg(Q)=2,\,ab\neq0$.
\end{rem}

\section{Quadratic polynomial automorphisms of $\C^3$}
Let $f$ be a quadratic polynomial automorphism of $\C^3$. Using
the classification of Forn\ae ss and Wu [FW], we show that -up to
conjugacy- $f$ or $f^2$ (or $f^{-1}$) is weakly regular. Moreover
$I^+$ (resp. $I^-$) is $f^{-1}$-attracting (resp. $f$-attracting)
except for certain mappings in the classes 4 or 5 below. Note that
$\l_1(f^{-1})=\l_2(f)$ since we are working in $\C^3$. Here
$\lambda_1(f)$ is the first dynamical degree of $f$,
$\lambda_1(f)=\lim[\delta_1(f^n)]^{1/n}$, where $\delta_1(f^n)$ is
the first algebraic degree of $f^n$ (see \cite{S99}).

\begin{thr}\label{t:class}
Let $f$ be a quadratic polynomial automorphism of $\C^3$ with
$\lambda_1(f) \neq \lambda_1(f^{-1})$. Then one of the following
holds:
\par 1) $f$ is conjugate to a regular automorphism with $X^-$
reduced to a point. In this case
$\lambda_1(f)=2<4=\lambda_1(f^{-1})$ and $I^{-}$ is
$f$-attracting.
\par 2) $f^2$ or $f^{-2}$ is conjugate to a mapping
from 1).
\par 3) $f$ is conjugate to
$$f(x,y,z)=(y[\alpha x+\beta y]+cx+dy+az,y^2+x,y)$$
where $a\alpha \neq 0$. In this case $f^{-1}$ is weakly regular
with $X^-=[0:0:1:0]$, $\lambda_1(f^{-1})=3 >2=\lambda_1(f)$, and
$I^-$ is $f$-attracting.
\par 4) $f$ or $f^{-1}$ is conjugate to
$$g(x,y,z)=(x^2-xz+c+y,az,bx+c'),$$
with $ab \neq 0$. In this case $g$ is weakly regular with
$X^+=[1:0:0:0]$, $\l_1(g)=2>\l_1(g^{-1})=(1+\sqrt{5})/2$, and
$I^+$ is $g^{-1}$-attracting if and only if $|b|<1$.
\par 5) $f$ is conjugate to
$$f(x,y,z)=(x[y+\alpha x]+az+c,x^2+dx+c'+by,x)$$
where $ab \neq 0$. In this case $f^{-1}$ is weakly regular,
$X^-=[0:0:1:0]$, $\lambda_1(f^{-1})=3 >2=\lambda_1(f)$, and $I^-$
is $f$-attracting if and only if $|b|>1$.
\end{thr}

\begin{pf} The quadratic polynomial automorphisms of
${\Bbb C}^3$ are classified into seven classes, up to affine
conjugacy \cite{FW98}. The growth of the degree of their forward
iterates is studied in \cite{BF00}. Two classes consist of affine
and elementary automorphisms $f$, so
$\lambda_1(f)=\lambda_1(f^{-1})=1$. We consider the remaining five
classes $H_1,\dots,H_5$ \cite{FW98}.

\par {\bf The classes $H_1$ and $H_2$.} By considering the degrees of
forward and backward iterates of the maps $H$ in these classes, it
is easy to see that $\lambda_1(H)=\lambda_1(H^{-1})\in\{1,2\}$.

\par {\bf The class $H_3$.} This class contains maps $H$ of the form
$$H(x,y,z)=(P(x,z)+a'y,Q(x)+z,x),\;\max\{deg(P),deg(Q)\}=2,\;a'\neq0.$$
We let $h=F\circ H\circ F^{-1}$, where $F(x,y,z)=(x,y-Q(z),z)$.
Then
\begin{equation}\label{special}
h(x,y,z)=(\alpha x^2+\alpha'xz+\alpha''z^2+c_1x+c_2z+c_3+a'y,z,x).
\end{equation}
The inverse map is
$$h^{-1}(x,y,z)=
\left(z,\frac{1}{a'}\,\left(x-\alpha
z^2-\alpha'yz-\alpha''y^2-c_1z-c_2y-c_3\right),y\right).$$ Using
the change of variables $(x,y,z)\ra(y,x,z)$ we see that $h^{-1}$
is conjugated to $h$, and the role of the coefficients
$\alpha,\,\alpha''$ interchanges. We have the following cases:

\par {\em Case A.} $\alpha\neq0\neq\alpha''$. Then
$\deg(h^n)=\deg(h^{-n})=2^n$, so
$\lambda_1(h)=\lambda_1(h^{-1})=2$.

\par {\em Case B.} $\alpha\neq0,\,\alpha''=0,\,\alpha'\neq0$. Then as
before $\deg(h^n)=2^n$ and $\lambda_1(h)=2$. The degrees of the
backward iterates $d_n=\deg(h^{-n})$ are given by Fibonacci's
numbers, $d_{n+2}=d_{n+1}+d_n$. So
$\lambda_1(h^{-1})=(1+\sqrt{5})/2$. Using the change of variables
$$F(x,y,z)=(\alpha x+v,\alpha a'y+s,-\alpha'z+r),\;
v=c_2\alpha/\alpha',\;r=2v-c_1,\;s=-\alpha a'r/\alpha',$$ we see
that $F\circ h\circ F^{-1}=g$, the map from $4)$. We have
$I^+(g)=\{t=x=0\}\cup\{t=x-z=0\}$ and $g(\{t=0\}\setminus
I^+)=X^+=[1:0:0:0]$. If $c=c'=0$ and $a=b^2$ the line
$\tau(\zeta)=(\zeta,b\zeta,\zeta)$ is $g$-invariant and
$g(\tau(\zeta))=\tau(b\zeta)$. So in this case $I^+$ is not
$g^{-1}$-attracting if $|b|\geq1$. We show in Lemma \ref{l:g4}
following this proof that $I^+$ is always $g^{-1}$-attracting if
$|b|<1$.

\par {\em Case C.} $\alpha\neq0,\,\alpha''=\alpha'=0$. Then $h^2$ is
regular, $\lambda_1(h^2)=4$, $\lambda_1(h^{-2})=2$, and
$X^+=[1:0:0:0]$.

\par {\em Case D.} $\alpha''\neq0,\,\alpha=0,\,\alpha'\neq0$. This is
similar to Case B, with the roles of $h$ and $h^{-1}$
interchanged, $\lambda_1(h)=(1+\sqrt{5})/2$ and
$\lambda_1(h^{-1})=2$.

\par {\em Case E.} $\alpha''\neq0,\,\alpha=\alpha'=0$. As in Case C,
$h^2$ is regular, $\lambda_1(h^2)=2$, $\lambda_1(h^{-2})=4$, and
$X^-=[0:1:0:0]$. The fact that $I^-$ is attracting for $f$ holds
for any regular automorphism $f$.

\par {\em Case F.} $\alpha=\alpha''=0,\,\alpha'\neq0$. As in Cases B
and D, $\lambda_1(h)=\lambda_1(h^{-1})=(1+\sqrt{5})/2$.

\par {\em Case G.} $\alpha=\alpha''=\alpha'=0$. Then $h$ is linear,
$\lambda_1(h)=\lambda_1(h^{-1})=1$.

\par {\bf The class $H_4$}. The maps $H$ in this class have the form
\begin{eqnarray*}
 & & H(x,y,z)=(P(x,y)+az,Q(y)+x,y),\;\max\{deg(P),deg(Q)\}=2,\;a\neq0,\\
 & & H^{-1}(x,y,z)=\left(y-Q(z),z,\frac{x}{a}+\widetilde P(y,z)\right),\;
\widetilde P(y,z)=-\frac{1}{a}\,P(y-Q(z),z).
\end{eqnarray*}
We write $P(x,y)=c_1x^2+c_2xy+c_3y^2+l.d.t.$,
$Q(y)=c_4y^2+l.d.t.$.

\par {\em Case A.} $c_4\neq0\neq c_1$. $H$ is regular,
$\lambda_1(H)=2,\,\lambda_1(H^{-1})=4,\, X^-=[0:0:1:0]$.

\par {\em Case B.} $c_4\neq0,\,c_1=0,\,c_2\neq0$. Then $H$ is
conjugated to the map $f$ of $3)$, $\lambda_1(f)=2$,
$\lambda_1(f^{-1})=3$, $f^{-1}$ is weakly regular,
$X^-=[0:0:1:0]$, $I^-$ is $f$-attracting (see \cite{C02}).

\par {\em Case C.} $c_4\neq0,\,c_1=c_2=0$. By \cite{CF99} (p.446)
either $H^2$ is regular, $\lambda_1(H^2)=4$,
$\lambda_1(H^{-2})=2$, $X^+=[c_3:c_4:0:0]$, or we have
$\deg(H^{\pm n})=2^n$.

\par {\em Case D.} $c_4=0$. If $F(x,y,z)=(x+Q(y),z,y)$, $F\circ H\circ
F^{-1}$ is the map from (\ref{special}).

\par {\bf The class $H_5$.} The maps $H$ in this class have form
\begin{eqnarray*}
 & & H(x,y,z)=(P(x,y)+az,Q(x)+by,x),\;\max\{deg(P),deg(Q)\}=2,\;a\neq0\neq b,\\
 & & H^{-1}(x,y,z)=\left(z,\frac{y-Q(z)}{b}\,,\frac{x}{a}+\widetilde P(y,z)\right),
\;\widetilde
P(y,z)=-\frac{1}{a}\,P\left(z,\frac{y-Q(z)}{b}\right).
\end{eqnarray*}
Let $P(x,y)=c_1x^2+c_2xy+c_3y^2+d_1x+d_2y+d_3$,
$Q(x)=c_4x^2+e_1x+e_2$.

\par {\em Case A.} $c_4\neq0\neq c_3$. $H$ is regular,
$\lambda_1(H)=2$, $\lambda_1(H^{-1})=4$, $X^-=[0:0:1:0]$.

\par {\em Case B.} $c_4\neq0,\,c_3=0,\,c_2\neq0$. Then $\deg(H^n)=2^n$
and $\deg(H^{-n})=3^n$. If
$$F(x,y,z)=(px+q,c_2y+r,pz+q),\;p^2=c_2c_4,\;q=pd_2/c_2, \;r=d_1-2qc_1/p,$$
then $F\circ H\circ F^{-1}$ is the map $f$ from $5)$,
$I^-=\{t=z=0\}$, $f^{-1}(\{t=0\}\setminus I^-)=X^-=[0:0:1:0]$. If
$|b|>1$ it is shown in \cite{GS} that $I^-$ is $f$-attracting. If
$|b|\leq1$ and if $f$ fixes the origin, then $f(0,y,0)=(0,by,0)$,
so $I^-$ is not $f$-attracting.

\par {\em Case C.} $c_4\neq0,\,c_3=c_2=0$. The inverse map is
$$H^{-1}(x,y,z)=\left(z,\frac{y-c_4z^2-e_1z-e_2}{b},\frac{x}{a}+
\frac{\gamma z^2}{a}-\frac{d_2y}{ab}+L(z)\right),$$ where
$\gamma=(d_2c_4/b)-c_1$ and $\deg(L)\leq1$. If
$c_1\neq0\neq\gamma$ then $\lambda_1(H)=\lambda_1(H^{-1})=2$. If
$c_1\neq0$ and $\gamma=0$ then $d_2\neq0$ and $H^2$ is regular,
$\lambda_1(H^2)=4$, $\lambda_1(H^{-2})=2$. If $c_1=0$ and
$d_2\neq0$ then $H^2$ is regular, $\lambda_1(H^2)=2$,
$\lambda_1(H^{-2})=4$. If $c_1=d_2=0$ then the degrees of all
iterates are bounded by 2.

\par {\em Case D.} $c_4=e_1=0$. If $c_1\neq0$ then
$\lambda_1(H)=\lambda_1(H^{-1})=2$. If $c_1=0$ then $deg(H^{\pm
n})\leq n+1$, so $\lambda_1(H)=\lambda_1(H^{-1})=1$.

\par {\em Case E.} $c_4=0,\,e_1\neq0$. We have that $F\circ H\circ
F^{-1}$ is the map $h$ from (\ref{special}), where
$$F(x,y,z)=\left(e_1x+by+e_2+\frac{e_2}{b}\,,-\frac{e_1z}{b}+
\frac{y}{b}\,,y+\frac{e_2}{b}\right).$$
\end{pf}

\begin{lem}\label{l:g4}
If $g(x,y,z)=(x^2-xz+c+y,az,bx+c')$ is the map from Theorem
\ref{t:class}, case $4)$, and $|b|<1$, then $I^+$ is
$g^{-1}$-attracting.
\end{lem}

\begin{pf} The inverse of $g$ has the form
$$g^{-1}(x,y,z)=(x_1,y_1,z_1)=\left(\frac{z}{b}+c'',\frac{z}{b}\left(
     \frac{y}{a}-\frac{z}{b}\right)+L(y,z)+x,\frac{y}{a}\right),$$
where $c''\in{\Bbb C}$ and $\deg(L)\leq1$. Recall that
$I^+=\{t=x=0\}\cup\{t=x-z=0\}$. We let $\a=|b|/(4|a|)$ and define
for $R>1$
\begin{eqnarray*}
 & & V_R=\left\{(x,y,z)\in\C^3:\,
 \max\{2\a|y|,|z|\}>\max\{2R,R^{1/3}|x|\}\right\},\\
 & & W_R=\left\{(x,y,z)\in\C^3:\,
 \max\{\a|y|,|x|\}>\max\{R,R^{1/3}|x-z|\}\right\}.
\end{eqnarray*}
Since $|b|<1$ we can find $\e>0$ such that $|b|<(1-2\e)/(1+\e)$.
The lemma follows if we show that for all $R$ sufficiently large
we have \begin{equation}\label{e:incl} g^{-1}(V_R)\subset
V_{2R}\cup W_{2R}\;, \;g^{-1}(W_R)\subset V_{2R}\cup W_{(1+\e)R}.
\end{equation}

\par We denote in the sequel by $C_g$ all constants which depend only
on the coefficients of $g$. For the first inclusion of
(\ref{e:incl}), let $(x,y,z)\in V_R$. We have two cases:

\par {\em Case A.} $2\a|y|\geq|z|$, so
$|y|>R/\a,\;|y|>R^{1/3}|x|/(2\a)$. We show that in this case
$g^{-1}(x,y,z)\in V_{2R}$. If $|y|/|a|>4R^{1/3}|z|/|b|$ then
$$2R^{1/3}|x_1|\leq2R^{1/3}\frac{|z|}{|b|}+2|c''|R^{1/3}<|z_1|\;,
\;|z_1|>\frac{R}{\a|a|}>4R.$$ If $|y|/|a|\leq4R^{1/3}|z|/|b|$,
using $|z|/|b|\leq2\a|y|/|b|=|y|/(2|a|)$, we get
$$|y_1|\geq\frac{|z|}{|b|}\left(\frac{|y|}{|a|}-\frac{|z|}{|b|}\right)
-|x|-|L(y,z)|\geq\frac{C_g|y|^2}{R^{1/3}}>\max\{4R,2R^{1/3}|x_1|\}.$$

\par {\em Case B.} $2\a|y|<|z|$, so $|z|>2R,\;|z|>R^{1/3}|x|$. If
$|x_1|>2R^{1/3}|x_1-z_1|$ then $g^{-1}(x,y,z)\in W_{2R}$, since
$|x_1|\geq|z|/|b|-|c''|>2R$. If $|x_1|\leq2R^{1/3}|x_1-z_1|$ then
$|z/b-y/a|\geq C_g|z|/R^{1/3}$, so $|y_1|>C_g|z|^2/R^{1/3}$ and
$g^{-1}(x,y,z)\in V_{2R}$.

\par To prove the second inclusion of (\ref{e:incl}), let $(x,y,z)\in
W_R$ and consider two cases:

\par {\em Case A.} $\a|y|\geq|x|$, so $|y|>R/\a,\;|y|>R^{1/3}|x-z|/\a$.
If $|z_1|>2R^{1/3}|x_1|$ then $g^{-1}(x,y,z)\in V_{2R}$, since
also $|z_1|=|y|/|a|>4R$. If $|z_1|\leq2R^{1/3}|x_1|$ then
$$\frac{|z|}{|b|}\geq\frac{|y|}{2|a|R^{1/3}}-|c''|\geq\frac{|y|}{3|a|R^{1/3}}
\;,\;\;|z|\leq|z-x|+|x|\leq\frac{\a|y|}{R^{1/3}}+\a|y|<2\a|y|.$$
It follows that $g^{-1}(x,y,z)\in V_{2R}$, since
$$|y_1|\geq
\frac{|y|}{3|a|R^{1/3}}\left(\frac{|y|}{|a|}-\frac{2\a|y|}{|b|}\right)-C_g|y|
>\frac{C_g|y|^2}{R^{1/3}}\;.$$

\par {\em Case B.} $\a|y|<|x|$, so $|x|>R,\;|x|>R^{1/3}|x-z|$. There
exists a large constant $M$ depending only on $g$, such that if
$|z/b-y/a|\geq M$ then $g^{-1}(x,y,z)\in W_{2R}$. Indeed, if $R$
is large we have $||z|-|x||<|x|/100$, so
$$\a|y_1|>\frac{|x|}{5|a|}\left|\frac{y}{a}-\frac{z}{b}\right|-C_g|x|\geq
\frac{|x|}{6|a|}\left|\frac{y}{a}-\frac{z}{b}\right|,$$ provided
that $M=M_g$ is sufficiently large. Therefore
$$\a|y_1|>\frac{RM}{6|a|}\geq2R\;,\;\;(2R)^{1/3}|x_1-z_1|\leq
2R^{1/3}\left|\frac{y}{a}-\frac{z}{b}\right|<\a|y_1|,$$ so
$g^{-1}(x,y,z)\in W_{2R}$. Finally, we assume that $|z/b-y/a|<M$.
For $R$ large we have $||z|-|x||<\e|x|$, so
$|x_1|\geq|z|/|b|-|c''|> (1-2\e)|x|/|b|>(1+\e)|x|$. Since $|x|>R$
and $|x_1-z_1|\leq M+|c''|$, we conclude that in this case
$g^{-1}(x,y,z)\in W_{(1+\e)R}$. \end{pf}

\end{document}